\documentclass[12pt]{article}
\usepackage{amsmath,amsxtra,amssymb,latexsym,amscd,amsthm}

\DeclareMathOperator{\wideg}{wideg}
\DeclareMathOperator{\Lwideg}{Lwideg}
\DeclareMathOperator{\Reso}{Res^o}
\DeclareMathOperator{\Resc}{Res^c}
\DeclareMathOperator{\Co}{\mathscr{C}^o}
\DeclareMathOperator{\Cc}{\mathscr{C}^c}
\DeclareMathOperator{\Res}{Res_x}
\DeclareMathOperator{\Syl}{Syl}
\DeclareMathOperator{\Image}{Im}
\DeclareMathOperator{\diag}{diag}
\DeclareMathOperator{\Ker}{Ker}
\DeclareMathOperator{\lcm}{lcm}

\usepackage[a4paper,left=2cm,right=2cm,top=2.0cm,bottom=2.0cm]{geometry}

\usepackage{tikz}
\usepackage{tikz-cd} 
\usepackage{array}
\usepackage{dsfont}
\usepackage{enumitem}
\usepackage{indentfirst}
\newtheorem{theorem}{Theorem}[section]
\newtheorem{lemma}{Lemma}[section]
\newtheorem{proposition}{Proposition}[section]
\newtheorem{corollary}{Corollary}[section]
\newtheorem{construction}{Construction}[section]
\newtheorem{algorithm}{Algorithm}[section]

\theoremstyle{definition}

\newtheorem*{remark}{Remark}

\usepackage{mathrsfs}
\begin{document}
	\title{ Applications of resultant of two $p$-adic power series }
	\author{Tran Hoang Anh\\ {Supervisor:} Dr. Zábrádi Gergely}
	\maketitle
\begin{abstract}
Given a prime $p$, and $v_p(a)$ stand for the $p$-adic valuation of the element $a$ in a finite extension $K$ of $\mathbf{Q}_p$, or more generally the field $\mathbf{C}_p$ which is the complete field of the algebraic closure $\mathbf{Q}_p$ with respect to the $p$-adic absolute value, denoted by $\lvert \cdot \rvert_p$. Let $F$ and $G$ be two ($p$-adic) power series with no common roots. We aim to estimate the maximal value $S$ and the minimal value $s$ of the function $\phi(x)=\min(v_p(F(x)), v_p(G(x)))$ over various domains, namely open and closed unit discs of $K$ or $\mathbf{C}_p$. To do this, we use partial resultants of two power series over certain domains defined by varied versions of the Weierstrass preparation theorem. Furthermore, the resultant of power series provides an efficient tool while studying the irreducibility and calculating the maximal value of $\phi$.
\end{abstract}

\section{Introduction and notations} 
Our object is in the range of the function $\phi$ which is thoroughly investigated in the case $F$ and $G$ are both polynomials in the paper $\cite{FZ}$. However, it gets tougher in terms of power series because the resultant is hardly defined. In this paper, we implement that kind of definition of power series $\mathcal{O}_K$ denote the ring of integers of $K$, and let $\pi$ be a uniformizer of $K$, then $\mathfrak{m}_K=\pi\mathcal{O}_K$ is the maximal ideal of $\mathcal{O}_K$. In the same manner, we define $\mathcal{O}_{\mathbf{C}_p}$ and $\mathfrak{m_{\mathbf{C}_p}}$ in $\mathbf{C}_p$ those are closed and open unit discs respectively.

Let $F(X)=\sum_{i=0}^{\infty}F_iX^i \in \mathbf{C}_p[\left[ X ]\right]$ converge in a region of $\mathbf{C}_p$, so-called the region of convergence of radius $\delta$. The Weierstrass preparation theorem (see $\cite{LB}$) claims the connection between the number of roots of $F(X)$ and the Weierstrass degree $\wideg(F)$ which is the smallest index $n$ such that  $F_n \in \mathcal{O}_{\mathbf{C}_p}^\times$, but in this paper, we extend it to be the smallest integer such that $v_p(F_n)= \min\{v_p(F_i) |\ i=0,1,2,\dotsc\}$, thus, the theorem remains unchanged as the following.

\begin{theorem} Let $F(X)$ be a power series of coefficients in $\mathcal{O}_K$ (or in $\mathcal{O}_{\mathbf{C}_p}$ generally), then $F(X)$ admits the unique form $F_{\wideg(F)}P(X)U(X)$ where $U(X) \in R_n[\left[ X ]\right]^{\times}$ and $P(X)= X^n+P_{n-1}X^{n-1}+\dotsc +P_0 \in R_n[X]$ is a distinguished polynomial for the ideal $I_n$ generated by $F_0,\dotsc,F_{n-1}$ of $R_n=\textbf{Z}\left[F_n,F_n^{-1},\{F_k\}_{k \geq n+1}\right][\left[F_0,\dotsc, F_{n-1}]\right]$.
	
\end{theorem}
Theorem 1.1 shows how the distinguished polynomial is able to replace the power series while studying the $p$-adic valuation on the open unit disc $\mathfrak{m_{\mathbf{C}_p}}$. For those converging on the closed unit disc, a bound for the number of roots in $\mathcal{O}_{\mathbf{C}_p}$ is given by the Strassman's theorem which also have a generalization known as the $p$-adic Weierstrass preparation theorem $\cite{FG}$ that again proves the existence of a polynomial such that for all $x \in \mathcal{O}_{\mathbf{C}_p}$, the polynomial and the original power series are of the same $p$-adic valuation.
\begin{theorem}
	Let $F(X)=\sum_{i=0}^{\infty}F_iX^i \in \mathbf{C}_p\left[\left[X\right]\right]$ be a power series such that $F(X)$ converges for $x \in \mathcal{O}_{\mathbf{C}_p}$ and $\Lwideg(F)$ be the integer defined by the following conditions\\
	$$ v_p(F_{\Lwideg(F)})=\min\{v_p(F_n)\}\  \text{and} \ v_p(F_n) > v_p(F_{\Lwideg(F)}) \ \forall n > \Lwideg(F)$$
	Then there exists a monic polynomial $g(X)$ of degree $F_{\Lwideg(F)}$ and a power series $h(X) \in \mathcal{O}_{\mathbf{C}_p}[\left[X\right]]^{\times}$, satisfying:
	\begin{itemize}
		\item[i)] $F(X)=F_{\Lwideg(F)}g(X)h(X)$,
		\item[ii)] The constant coefficient of $h(X)$ is a unit and the other coefficients are not,
		\item[iii)] $h(X)$ has no zero on the closed unit disc,
		\item[iv)] If coefficients of $F(X)$ belongs to a finite extension $K$, then so are all coefficients of $g(X)$ and $h(X)$.
	\end{itemize}
	
\end{theorem} 
The two versions of Weierstrass preparation theorem provide an essential tool to construct the resultant of two power series (see $\cite{LB}$) that eventually raises various bounds of $S$. In the present paper, we approach the function $\phi$ by two different methods, one of them is a refined study of the paper $\cite{FZ}$. The polynomials $P$ in theorem 1.1 and $g$ in theorem 1.2 are so-called distinguished polynomials regarding the open and closed unit discs respectively.
\section{Resultants of $p$-adic power series} 
The theory of Newton polygons declares that both distinguished polynomials defined above have all roots in $\mathcal{O}_{\mathbf{C}_p}$ and $\mathfrak{m}_{\mathbf{C}_p}$ respectively in the cases of open and closed discs. Furthermore, both the distinguished polynomials cover all roots of the power series on discs. Therefore, we can consider resultants in usual ways. Because of the similarity of two cases open and closed discs, we shall treat those in parallel.
\begin{theorem}
	Let $F(X)$ and $G(X)$ be power series in $\mathbf{C}_p[\left[X\right]]$ and converge in $\mathfrak{m}_{\mathbf{C}_p}$. In addition, $\wideg(F)=n$ and $\wideg(G)=N$ are finite, then there exists 
	$$\Reso(F,G)=F_n^{N}\prod_{\substack{z\in \mathfrak{m_{\mathbf{C}_p}} \\F(z)=0}}G(z)$$
	In particular, $\Reso(F,G) \in \textbf{Z}[F_n,F_n^{-1}, \{F_k\}_{k>n},\{G_k\}][\left[F_0,\dotsc,F_{n-1}\right]]$ if $F_n$ is a unit.\\
	In the case the regions of convergence contains $\mathcal{O}_{\mathbf{C}_p}$, then $\Lwideg(F)=m$ and $\Lwideg(G)=M$ are finite and there exists 
	$$\Resc(F,G)=F_m^M\prod_{\substack{z \in \mathcal{O}_{\mathbf{C}_p}\\F(z)=0}}G(z)$$
	The condition  $\Resc(F,G) \in \textbf{Z}[F_m,F_m^{-1}, \{F_k\}_{k>m},\{G_k\}][\left[F_0,\dotsc,F_{m-1}\right]]$ also holds if $F_m$ is a unit.
\end{theorem}
One should mention here is that the region of convergence is either open and closed, thus, the corresponding resultants should be defined in varied manners. The proof of the first part was represented in [Theorem 7.3, \cite{LW}] by using Washington's theorem for the case $F$ and $G$ of coefficients in a finite extension of $\mathbf{Q}_p$. Following the method with several updates also provides a solution for the whole theorem. 
\begin{theorem}[\textbf{A generalization of Washington's theorem}]
	Let $f$ and $g$ be power series with coefficients in $\mathcal{O}\left[\left[X\right]\right]$ (the ring $\mathcal{O}$ here can be $\mathcal{O}_K$ or $\mathcal{O}_{\mathbf{C}_p}$), and $f(X)=\sum_{i=0}^{\infty}a_iX^i$.
	\begin{itemize}
		\item[1)] If $f$ converges in $\mathfrak{m_{\mathbf{C}_p}}$, $\wideg(f)=n$ and $a_n$ is a unit then there exists a unique form $g=fq+r$ where $r \in \mathcal{O}\left[X\right]$ of degree at most $n-1$ and $q \in \mathcal{O}\left[\left[X\right]\right]$.
		\item[2)] If $f$ and $g$ converge in $\mathcal{O}_{\mathbf{C}_p}$, $\Lwideg(f)=m$ and $a_m$ is a unit then there exists a unique form $g=fq+r$ where $r \in \mathcal{O}\left[X\right]$ of degree at most $m-1$ and $q \in \mathcal{O}\left[\left[X\right]\right]$.
	\end{itemize}
\end{theorem} 
\begin{proof}  We separate into two cases: 
	\begin{itemize}
		\item[1)] The uniqueness is easily shown by indirectly assuming $fq+r=0$ and $q,r \neq 0$. Without loss of generality, let $t=\min\{v_p(a_i), i=1,\dotsc,n-1\}$ then either $p^t \nmid q$ or $p^t \nmid r$ happens. Because the degree of $r$ is at most $n-1$ then $p^t \nmid qf$. Since $f$ has the coefficient $a_n$ is a unit, so $p^t \mid q$ that implies a contradiction.\\
		The existence is proved by a "shift operator" $\tau_n$ over the power ring $\mathcal{O}$$\left[\left[X\right]\right]$ defined by
		$$\tau_n\left(\sum_{i=0}^{\infty}b_iX^i\right)=\sum_{i=n}^{\infty}b_iX^{i-n}$$
	    Let $U(X)=\tau_n(f)$ and $P(X)=f(X)-x^nU(X)$ then $U(X)$ is a unit power series and  coefficients of $P$ is divisible by $p^t$. We define the power series $q$ by 
	    $$q=U^{-1}\sum_{i=0}^{\infty}(-1)^i\left(\tau_n\circ PU^{-1}\right)^{i}(\tau_n(g))$$
	    This power series is well-defined by $i$-th term of sum is divisible by $p^{ti}$. Therefore, one can be obtained is $\tau_n(g-qf)=0$ that makes $r=g-qf$ is a polynomial of degree at most $n-1$.
	    \item[2)] There are no such polynomials $q$ and $r$ that degree of $r$ is at most $m-1$ and $qf+r=0$. If there exist, then comparing the coefficients of degree at least $m$ shows that every coefficient of $q$ is not a unit, thus are coefficients of $r$. But it is a contradiction with the initial assumption of existence of a unit among coefficients of either $q$ or $r$.\\
	    The existences of $q$ and $r$ are also the same as the above case by replacing $n$ with $m$. The only problem here is the convergence of coefficients of $q$. Consider the power series $f$ converging in $\mathcal{O}_{\mathbf{C}_p}$, so $v_p(a_i) \to \infty$ as $i \to  \infty$. Let $g(X)=\sum_{i=0}^{\infty}b_iX^i$ and $U^{-1}(X)=\sum_{i=0}^{\infty}c_iX^i$, thus, $U^{-1}$ has the same region of convergence as $U$ which also implies $v_p(c_i) \to \infty$ as $i \to  \infty$. For any natural number $l$, the coefficient of $\left(\tau_m\circ PU^{-1}\right)^{k}(\tau_m(g))$ is 
	    $$\sum b_i\prod_{j=1}^{k}a_{v_j}c_{u_j}$$
	    where $i+\sum_{j=1}^{k}v_j+u_j=l+km$, $0 \leq v_j \leq m-1$ and $i+\sum_{j=1}^{w}v_j+u_j \leq (w+1)m \ \forall w=\overline{0,k}$.
	    Since $\Lwideg(f)=m$ then every coefficient of $U$ but the constant is not unit and $U$ also converges in $\mathcal{O}_{\mathbf{C}_p}$, so let $\min\{v_p(a_i) \ | \ i\geq m\}=t >0$. It is prominent that all coefficients of $U^{-1}$ except the constant one is also divisible by $p^t$. Combining with the $p$-adic valuation of coefficient of $g$ tending to $+\infty$, one is obtained $v_p(\sum b_i\prod_{j=1}^{k}a_{v_j}c_{u_j}) \to +\infty$ as $k \to +\infty$ which ensures the existence of $q$.
	\end{itemize}	 
	
\end{proof}
The theorem gives direct proofs for theorem 1.1 and 1.2 by taking $P=\dfrac{F-xFq}{c \cdot F_{\wideg(F)}}$ where $c$ is the leading coefficient of $\dfrac{F-xFq}{ F_{\wideg(F)}}$ or $P=2X^{\wideg(F)}-q'\cdot \dfrac{F}{F_{\wideg(F)}}$ alternatively where $q'$ is the division power series of $X^{\wideg(F)}$ applying theorem 2.2. Replacing $\wideg$ by $\Lwideg$, the distinguished polynomial in theorem 1.2 is defined in the same manner. Let $\Co(F)$ and $\Cc(F)$ be the distinguished polynomials corresponding to the open and closed unit discs.

We consider the $p$-adic valuation of the resultants of power series. As the definition, it is the product of the power series $G$ at zeros of $P$ on the (open or closed discs), but one is shown that the $p$-adic valuation depends on the distinguished polynomials as the following.
\begin{proposition}
	Let $F$ and $G$ be power series in $\mathcal{O}_{\mathbf{C}_p}\left[\left[X\right]\right]$ converging on the unit disc (either open or closed) and have distinguished polynomials $P$ and $Q$ respectively, then the $p$-adic valuation of the resultant is the same as $\Res(F_nP,G_NQ)$ where $n$ and $N$ are $\wideg$ or $\Lwideg$ with respect to $F$ and $G$ respectively.\\
	In particular, $$ \Res(F_nP,G_NQ)=F_n^NG_N^n\prod_{\substack{1\leq i\leq n\\1\leq j \leq N}} (z_i-w_j)=\det\left(\Syl(F_nP,G_NQ,x)\right)$$
	where $\{z_i| i=\overline{1,n}\}$ and $\{w_j|j=\overline{1,N}\}$ are the sets of roots in the unit disc of $F$ and $G$, the matrix $\Syl(F_nP,G_NQ,x)$ is the Sylvester matrix of two polynomials. 
	Furthermore, the supermum of the function $\phi$ over the unit disc domains $S$ admits an upper bound.
	$$S \leq v_p(\Res(F_nP,G_NQ))$$
\end{proposition}
\begin{proof}
	The first part is a consequence of the above theorem and the inequality is proven in the same method if $P$ and $Q$ are polynomials with integer coefficients, which is represent in $\cite{PJ}$.
\end{proof}
As usual, the resultant of power series implies those have common root or not in the unit disc. Moreover, We possibly determine when two power series have common root in the intersection of regions of convergence, namely $B(0,\delta)$ or $\overline{B}(0,\delta)$, and let $v_p(\delta)=\mu$. One should be careful here is that $\delta$ is not necessarily in $\mathbf{C}_p$ according to its definition, and the number $\mu$ is defined by taking the limit regarding to $p$-adic valuation. In section 4, we shall work on it deeper, but now we pay attention to the property of $\mu$ which is a real number and for $x \in \mathbf{C}_p$ such that $v_p(x) > \mu$ then both power series $F$ and $G$ converge at the point $x$. Especially, if $\mu \in Q$,  a point $x$ of $p$-adic valuation $\mu$ could be a converging point depending on the initial condition of $F$ and $G$. Furthermore, we study the power series $F \circ p^\mu$ and $G \circ p^\mu$ which have the intersection of regions of convergence the unit disc(either open or closed). For $\mu \notin Q$, we consider $F \circ p^{y_n}$ and $G \circ p^{y_n}$ over the unit disc, and let $(y_n)$ be a rational number sequence that strictly decreases to $\mu$. Hence, $F$ and $G$ have a common root if and only if there is an integer $n$ such that $F \circ p^{y_n}$ and $G \circ p^{y_n}$  have common root in the unit disc which equivalently implies the resultant is zero.
\section{Power series with coefficients in $\mathcal{O}_K$}
In this section, our object is to give an upper bound on the maximal value $S(K)$ of the function $\phi$ on the domain $\mathcal{O}_K$ where $K$ is a finite extension of $\mathbf{Q}_p$ with ramification index $e$ and residue degree $f$, $\pi$ is a uniformizer. The method to do this follows the paper $\cite{FZ}$.
\begin{theorem}
	Let $F$ and $G$ be power series with coefficient in $\mathcal{O}_K$ of finite $\wideg$, namely $\wideg(F)=n$ and $\wideg(G)=m$. In addition, $F_n$ and $G_m$ are units in $\mathcal{O}_{\mathbf{C}_p}$ 
	\begin{itemize}
		\item[i)] The ring generated by all power series converging in $\mathfrak{m_{\mathbf{C}_p}}$ in $\mathcal{O}_K\left[\left[X\right]\right]$ is denoted by $\mathcal{O}_K^o\left[\left[X\right]\right]$ and the resultant on the open unit disc is nonzero, then 
		$$\left|\mathcal{O}_K^o\left[\left[X\right]\right]/(F,G)\right|=p^{f(v_{\pi}(\Reso(F,G)))}$$
		\item[ii)] The ring generated by all power series converging in $\mathcal{O}_{\mathbf{C}_p}$ in $\mathcal{O}_K\left[\left[X\right]\right]$ is denoted by $\mathcal{O}_K^c\left[\left[X\right]\right]$ and the resultant on the closed unit disc is nonzero, furthermore, if $\Lwideg(F)=N$ and $\Lwideg(G)=M$ then 
		$$\left|\mathcal{O}_K^c\left[\left[X\right]\right]/(F,G)\right|=p^{f(v_{\pi}(\Resc(F,G)))}$$
	\end{itemize} 
\end{theorem}
\begin{proof}The two parts are proved in the same method, here is the solution for part $i$. The first part of theorem 2.2 claims the division with the remainder on $\mathcal{O}_K^o\left[\left[X\right]\right]$, thus, every power series in $\mathcal{O}_K^o\left[\left[X\right]\right]$ admits the unique form $U+FGH$ where $H$ is a power series and $U$ is a polynomial of degree at most $\wideg(FG)-1=n+m-1$. Besides, theorem 1.1 implies that $(F,G)=(\mathscr{C}^o(F),\mathscr{C}^o(G))$.
	
Let construct a correspondence between the vector space $\mathcal{O}_K^{m+n}$ and the vector space $\mathcal{O}_K\left[X\right]_{< m+n}$ of polynomials of degree less than $m+n$ by mapping the standard basis with the basis $\{1,X,\dotsc,X^{n+m-1}\}$. Consequently, the subgroup $(F,G)=(\mathscr{C}^o(F),\mathscr{C}^o(G))$ is identified with the subgroup generated by the set of rows of the matrix $T$ which is the $\Syl(\mathscr{C}^o(F),\mathscr{C}^o(G))$ defined by 
		$$\begin{pmatrix}
			a_0&a_1&\dotsc&\dotsc&a_n& \ & \ & \ &\\
			\ &a_0&a_1&\dotsc&\dotsc&a_n& \ & \ & \\
			\ & \ &\dotsc&\dotsc&\dotsc&\dotsc&\dotsc& \ & \\
			\ & \ &\ \ &a_0&a_1&\dotsc&\dotsc&a_n & \\
			b_0&b_1&\dotsc&\dotsc&b_m& \ & \ & \ &\\
			\ &b_0&b_1&\dotsc&\dotsc&b_m& \ & \ & \\
			\ & \ &\dotsc&\dotsc&\dotsc&\dotsc&\dotsc& \ & \\
			\ & \ &\ \ &b_0&b_1&\dotsc&\dotsc&b_m & \\
			
		\end{pmatrix}$$
where $\mathscr{C}^o(F)=\sum_{i=0}^{n}a_iX^i$ and $\mathscr{C}^o(G)=\sum_{i=0}^{m}b_iX^i$\\
Hence, $$\left|\mathcal{O}_K^o\left[\left[X\right]\right]/(F,G)\right|=\mathcal{O}_K\left[X\right]_{< m+n}/(\mathscr{C}^o(F),\mathscr{C}^o(G))=\mathcal{O}_K^{m+n}/\langle\text{rows of T}\rangle=\mathcal{O}_K^{m+n}/\Image(T)$$
In other hand, the ring $\mathcal{O}_K$ is a principal ideal domain, therefore, the Smith normal form $\cite{KR}$ works on the matrix $T$, namely 
$$T=A \diag(f_1,\dotsc,f_{n+m}) B$$
where the matrices $A$ and $B$ are invertible, and for each $i=\overline{1,n+m-1}$, $f_{i-1} | f_i$. Thus, $$v_\pi(\Reso(F,G))=v_\pi(\Res(\mathscr{C}^o(F),\mathscr{C}^o(G)))=v_\pi(\det(T))=v_\pi\left(\prod_{i=1}^{n+m}f_i\right)$$ and $$\mathcal{O}_K^{m+n}/\Image(T)\simeq\mathcal{O}_K^{m+n}/\Image(\diag(f_1,\dotsc,f_{n+m}))\simeq\bigotimes_{i=1}^{n+m}\mathcal{O}_K/(f_i)$$
Furthermore, $$\left|\mathcal{O}_K/(f_i)\right|=p^{f\cdot v_\pi(f_i)}$$
Combining what we have obtained provides what we need.

\end{proof}
\begin{corollary}
\begin{itemize} Bounds of the function $\phi$ in both cases:
		\item[i)] Let $S(\mathbf{C}_p)$ be the supermum of the function $\phi$ on the domain $\mathfrak{m}_{\mathbf{C}_p}$ in the case $F$ and $G$ in  $\mathcal{O}_K^o\left[\left[X\right]\right]$, then $S(\mathbf{C}_p) \leq f_{n+m}$
		\item[ii)] Let $S(\mathbf{C}_p)$ be the supermum of the function $\phi$ on the domain $\mathcal{O}_{\mathbf{C}_p}$ in the case $F$ and $G$ in  $\mathcal{O}_K^c\left[\left[X\right]\right]$, then $S(\mathbf{C}_p) \leq f_{n+m}'$ which is defined similarly in the case of closed unit disc.
\end{itemize}
	 
\end{corollary}
\begin{proof} It is enough to prove the first one.
	We define the $p$-adic valuation on $\mathcal{O}_{\mathbf{C}_p}^{n+m}$ in the manner that for $x=(x_1,x_2,\dotsc,x_{n+m})$, $v_p(x)=\min\{v_p(x_i)\}$. One can be shown here is that for $x \in \mathfrak{m}_{\mathbf{C}_p}$, then $\phi(x)=v_p(T\cdot(x^{n+m-1},\dotsc,1)^T)$. Thus, $$S(\mathbf{C}_p) \leq \sup\{v_p(Tx)\ | \ x \in \mathcal{O}_{\mathbf{C}_p}^{n+m}, v_p(x)=0\}$$
	Notice that the matrices $A$ and $B$ are invertible, which means they preserve the $p$-adic valuation of vectors, therefore, $T$ can be replaced by $\diag(f_1,\dotsc,f_{n+m})$ which also implies that $S(\mathbf{C}_p) \leq f_{n+m}$.
\end{proof}
In the paper $\cite{FZ}$, there are many numerical functions that serve the estimations of the maximal value of the greatest common divisor, thus, we try to extend it to the finite extension $K$ of $\mathbf{Q}_p$. The work below is a generalization of the paper $\cite{KR}$ to counting polynomial functions mod $\pi^n$.

Let define $R_n$ be the ring of polynomial functions $\mathcal{O}_K/(\pi^n) \to \mathcal{O}_K/(\pi^n)$ and its cardinality is calculated by the functions defined on $\textbf{Z}$
$$ \alpha(j)= \sum_{k=1}^{\infty} \Bigl\lfloor\dfrac{j}{p^{kf}}\Bigr\rfloor, \ \text{$f$ here is the residue degree of the field $K$}$$ and $\beta(m)=\min\{j: \alpha(j) \geq m\}$. Put $B(s)=\sum_{m=1}^{s}\beta(m)$.

Before calculating the cardinality, we analogize the ring $\mathcal{O}_K/(\pi^n)$ with the ring $\textbf{Z}$ by giving it an order. Let quickly remind prominent properties of the field $K$ which is the order of residue field $k=\mathcal{O}_K/(\pi)$ is $p^f$, and let $\{c_0,c_1,\dotsc, c_{p^f-1}\}$ be a fixed set of representatives and $c_0$ represent the element zero, then any element $X$ in $K$ has a unique $\pi$-adic expansion form $X=\sum_{i=-m}^{\infty}c_i\pi^i$.Therefore, let a positive integer $j$ of the form $j=\sum_{i=0}^{f(n-1)}a_ip^{fi}$ in the $p_f$-base where $a_i\in\{0,1,\dotsc, p^f-1\}$ then denote the $j$- element of $\mathcal{O}_K/(\pi^n)$ by $d_j=\sum_{i=0}^{n-1}a_i\pi^{i}$.

Let $X^{[j]}=(X-d_0)(X-d_1)\dotsm(X-d_{j-1})$ for positive integer $j$ and define $x^{[0]}=d_1$. It is easily seen that $X^{[j]} \in R_n$ and the set $\{X^{[j]},X^{[1]},\dotsc,X^{[p^nf-1]}\}$ form a $\mathcal{O}_K$-basis of $R_n$.
\begin{lemma}
	For any positive integer $j$, $v_{\pi}\left(X^{[j]}\right) \geq \alpha(j) \ \forall X \in \mathcal{O}_K$ and $v_{\pi}\left(d_{j+1}^{[j]}\right)=\alpha(j)$.
\end{lemma}
\begin{proof}
	The ideal is the same as in the ring of integer numbers. 
	\begin{itemize}
		\item First, any $p^f$ consecutive element $d_i$'s forms a point-wise incongruent set there for each of them contributes $\Bigl\lfloor\dfrac{j}{p^{f}}\Bigr\rfloor$ to the $\pi$-adic valuation.
		\item The same happens for any  $p^{kf}$ consecutive element $d_i$'s.
	\end{itemize}
Besides, the $\pi$-adic valuation of $d_{j+1}^{[j]}=\alpha(j)$ is the same as $\prod_{i=1}^{j}d_i$ which is exactly $\alpha(j)$.
\end{proof}
\begin{theorem}
	For any polynomial $F \in R_n$, there exists uniquely a polynomial with $F=F_n$ and $F_n=\sum_{i+\alpha(j)<n}b_{ij}\pi^iX^{[j]}$ where $b_{ij} \in \{c_0,c_1,\dotsc, c_{p^f-1}\}$.
\end{theorem}
\begin{proof}
	The lemma 3.1 implies that if $i+\alpha(j) \geq n$ then the polynomial $\pi^ix^{[j]}$ vanishes mod $\pi^n$, combining with the $\mathcal{O}_K$-basis of $R_n$ which is $\{X^{[j]},X^{[1]},\dotsc,X^{[p^nf-1]}\}$ ensures the existence of $F_n$. 
	
	Indirectly assume that there is a non-zero polynomial of the form $\sum_{i+\alpha(j)<n}b_{ij}\pi^iX^{[j]}$ vanishing mod $\pi^n$. Rewriting $e_j=\sum_{i=0}^{n-1-\alpha(j)}b_{ij}\pi^i$ and replacing $X$ with $d_k$ with $k$ runs from 1 to $p^{nf}-1$ consequently claims that all $e_j$'s vanish mod $\pi^n$ which equivalently means that $b_{ij}=0 \ \forall i$. It is a contradiction with our hypothesis.
\end{proof}
Let define the map $\varphi_n: R_n \to R_n-1$ by $\varphi(F_n)=F_n-1$, then $\varphi_n$ is a surjective homomorphism from $R_n$ onto $R_{n-1}$.
\begin{theorem}
	Every element in $\Ker(\varphi_n)$ is uniquely of the form $\sum_{i+\alpha(j)=n-1}b_{ij}\pi^iX^{[j]}$ where $b_{ij} \in \{c_0,c_1,\dotsc, c_{p^f-1}\}$. Furthermore, the cardinality of $\Ker(\varphi_n)$ is precisely $p^{f\beta(n)}$
\end{theorem}  
\begin{proof}
	When we rewrite the polynomial $F_n$ by the form in $R_{n-1}$ by the above method, those terms $b_{ij}\pi^iX^{[j]}$ with $i+\alpha(j)<n-1$ remain unchanged. In contrast, one with $i+\alpha(j)=n-1$ vanishes mod $\pi^{n-1}$. Therefore, elements of $\Ker(\varphi_n)$ are of the form $\sum_{i+\alpha(j)=n-1}b_{ij}\pi^iX^{[j]}$.
	
	The equation $i+\alpha(j)=n-1$ of variable $i$ with fixed integer $0 \leq j \leq \beta(n)-1$  has exactly one solution, hence, the number of  $\Ker(\varphi_n)$ is $p^{f\beta(n)}$.
\end{proof}
\begin{corollary}
	The cardinality of $R_n$ is $p^{fB(n)}$.
\end{corollary}
\begin{proof}
	The result is easily obtained by induction.
\end{proof}
\begin{remark}
	For the need of estimations, we shall also need to study the cardinality of $R'_n$ which is the ring of polynomial functions $\mathfrak{m}_K/(\pi^n) \to \mathfrak{m}_K/(\pi^n)$. The above calculations still work with modifications as following
	$$ \alpha'(j)= j+\sum_{k=1}^{\infty} \Bigl\lfloor\dfrac{j}{p^{kf}}\Bigr\rfloor, \ \text{$f$ here is the residue degree of the field $K$}$$ and $\beta'(m)=\min\{j: \alpha'(j) \geq m\}$. Put $B'(s)=\sum_{m=1}^{s}\beta'(m)$. Then 
	$$\left|R'_n\right|=p^{fB'(n)}$$
	
\end{remark}
Let $S_K$ and $s_K$ be maximal and minimal values of the function $\phi$ over the domains  $\mathcal{O}_K$ respectively and $S=eS_K, s=es_K$ when replacing $p$-adic valuation by $\pi$-adic valuation. We use the notations $S'$ and $s'$ in the case of domain $\mathfrak{m}_K$. We observe the ideal $$I'_{S',s'}=\{F \in \mathcal{O}_K^o\left[\left[X\right]\right]: \pi^{s'} \mid F(X) \ \forall X \in \mathfrak{m}_K, \pi^{S'} \mid F(0) \}$$ of the ring $\mathcal{O}_K^o \left[\left[X\right]\right]$. Put $R'_{S',s'}=\mathcal{O}_K^o\left[\left[X\right]\right]/I'_{S',s'}$\\
and the ideal $$I_{S,s}=\{F \in \mathcal{O}_K^c\left[\left[X\right]\right]: \pi^s \mid F(X) \ \forall X \in \mathcal{O}_K, \pi^S \mid F(0) \}$$ of the ring $\mathcal{O}_K^c\left[\left[X\right]\right]$. Put $R_{S,s}=\mathcal{O}_K^c\left[\left[X\right]\right]/I_{S,s}$.
\begin{lemma}
	$$\left|R'_{S',s'}\right|=p^{f\left(S'-s'+B'(s)\right)}, \ \left|R_{S,s}\right|=p^{f\left(S-s+B'(s)\right)}$$ 
\end{lemma}
\begin{proof}
	In the case $S=s$ as well as $S'=s'$, $R'_{S',s'}$ and $R_{S,s}$ are $R'_{s'}$ and $R_s$ respectively and the equations have been proved.
	
	For the case $S>s$, by considering the map $I_{s,s} \to \mathcal{O}_K/(\pi^S)$ defined by $F \mapsto F(0)$, the image of this map is $(\pi^s)/(\pi^S)$ and kernel is $I_{S,s}$, whence $\left|I_{S,s}/I{s,s}\right|=p^{f(S-s)}$ which also leads to $\left|R_{S,s}\right|/\left|R_{s,s}\right|=p^{f(S-s)}$. It also works on $R'_{S',s'}$.
\end{proof}
The first main estimation is the refined study of the paper $\cite{FZ}$.

\begin{theorem}
	Let $F$ and $G$ be power series with coefficient in $\mathcal{O}_K$ of finite $\wideg$, namely $\wideg(F)=n$ and $\wideg(G)=m$. In addition, $F_n$ and $G_m$ are units in $\mathcal{O}_{\mathbf{C}_p}$.
	\begin{itemize}
		\item[i)] If $F$ and $G$ in $\mathcal{O}_K^o\left[\left[X\right]\right]$ then for all non-negative integers $t$ the following inequality holds
		$$v_{\pi}\left(\Reso(F,G)\right)-S' \geq B'(s'+t)-2B'(t)-s'$$ equivalently
		$$v_p\left(\Reso(F,G)\right)-S'_K \geq \dfrac{1}{e}\left(B'(es'_K+t)-2B'(t)-es'_K\right)$$
		\item[i)] If $F$ and $G$ in $\mathcal{O}_K^c\left[\left[X\right]\right]$ then for all non-negative integers $t$ the following inequality holds
		$$v_{\pi}\left(\Resc(F,G)\right)-S \geq B(s+t)-2B(t)-s$$ equivalently
		$$v_p\left(\Reso(F,G)\right)-S_K \geq \dfrac{1}{e}\left(B(es_K+t)-2B(t)-es_K\right)$$
	\end{itemize}

\end{theorem}
\begin{proof}
	As usual, we only need to prove one of them and the other is mostly the same. The one proved here is part $ii$. Because the resultant of power series is also translation-invariant as the normal one, we may assume that the maximal value achieved at the point zero. Notice that 
	$$\mathcal{O}_K^c\left[\left[X\right]\right]/(F,G) \supseteq \mathcal{O}_K^c\left[\left[X\right]\right]/\left((F,G)+I_{S+t,s+t}\right)=R_{S+t,s+t}/\left(\overline{F},\overline{G}\right)$$
	with $\overline{F}$ and $\overline{G}$ are the natural images in $R_{S+t,s+t}$ by $F$ and $G$, respectively which is annihilated by $I_{t,t}$, then 
	$$v_{p^f}\left(\left|\left(\overline{F},\overline{G}\right)\right|\right)=v_{p^f}\left(\left|\left(\overline{F}\right)\right|\right)+v_{p^f}\left(\left|\left(\overline{G}\right)\right|\right)-v_{p^f}\left(\left(\overline{F}\right) \cap \left(\overline{G}\right)\right)\leq 2B(t)$$ 
	Using theorem 3.1 and lemma 3.2, we obtain 
	$$v_{\pi}\left(\Resc(F,G)\right) \geq (S+t)-(s-t)+B(s+t)-2B(t)$$ as we desire.
	
\end{proof}
To show the sharpness of the inequality, we present a construction as following for a fixed $p$.
\begin{construction}
	Let $$F(X):= \prod_{j=0}^{\beta(s)-1}(X-d_j);$$ and $$G(X):=\pi^s+\prod_{i=0}^{p^f-1}(X-c_i)^{s+1}$$
	for a positive integer $s$. Then over the domain $\mathcal{O}_K$, the function $\phi \equiv s$. The $p$-adic valuation of resultant is $s\beta(s)$, whence $v_{\pi}\left(\Resc(F,G)\right)-s=s(\beta(s)-1)$ which is approximate to $(p_f-1)s^2$ as $s$ tends to $\infty$.
\end{construction}
\begin{proof}
	As $s \to \infty$ $B(s) \sim (p^f-1)\dfrac{s^2}{2}$, thus, if we let $s=t$ in the inequality, $B(2s)-2B(s)-s \sim (p^f-1)s^2$.
\end{proof}
\begin{theorem}
	Let $F$ and $G$ be power series of the same conditions as above. In addition, $s$ or $s'$ is not larger than $p^f-1$ then 
	\begin{itemize}
		\item[i)] If $F$ and $G$ in $\mathcal{O}_K^o\left[\left[X\right]\right]$ then 
		$$v_{\pi}\left(\Reso(F,G)\right)-S' \geq s'(s'-1)$$
		\item[ii)] If $F$ and $G$ in $\mathcal{O}_K^c\left[\left[X\right]\right]$ then 
		$$v_{\pi}\left(\Resc(F,G)\right)-S \geq p^fs^2-s$$
	\end{itemize}
\end{theorem}
\begin{proof}
	The proof is the same in the case of integer numbers which is represented in the paper $\cite{FZ}$.
\end{proof}
\begin{corollary}
	In the case $F_{\wideg(F)}$ and $G_{\wideg(G)}$ might not be unit, then we still use the above inequality to estimate the range of the function $\phi$ by dividing both of them by the coefficients of degree wideg. In particular, $$\phi(F,G)(X) \leq \max\{v_p(F_{\wideg(F)}),v_p(G_{\wideg(G)})\}+\phi\left(\dfrac{F}{F_{\wideg(F)}},\dfrac{G}{G_{\wideg(G)}}\right)(X)$$
\end{corollary}
 We have just estimated the maximal value $S$ of the function $\phi$ over finite extension of $\mathbf{Q}_p$. In fact, we can also estimate the minimal $s$ value of this function by the theorem
 \begin{theorem}
 	Let $P$ be a monic polynomial of degree $n$ and coefficients in $K$. Assume that $\pi^s \mid P(X) \ \forall X \in \mathcal{O}_K$, then $n \geq \beta(s)$.\\
 	Furthermore, If $P$ is a power series of the same property and $\Lwideg(P)$ is a unit then $\Lwideg(P) \geq \beta(s)$
 \end{theorem} 
\begin{proof}
	Assume indirectly that $n < \beta(s)$ then $P=P_s$ mod $\pi^s$ and have the coefficient of $x^{[n]}$ is $1$ in $R_s$, therefore consider $P(d_i)$ with $i$ run from $0$ to $n$, we have all coefficient of $x^{[i]}$ must be divisible by $\pi$, a contradiction. The rest is obtained directly from what we just proved.
\end{proof}
\section{Analyzing roots of $p$-adic power series}
In this section, we discuss deeper the problem we have left in section 2. The convergence radius $\delta$ of a power series $F(X) \in \mathbf{C}_p\left[\left[X\right]\right]$ is a real number such that for any $\left|X\right|_p < \delta$ then $F(X)$ converges. Moreover, $\delta=\left(\limsup_{n \to \infty} \sqrt[n]{\left|F_n\right|_p}\right)^{-1}$ which is of the form $p^\mu$ where $\mu=-\liminf_{n \to \infty}\dfrac{v_p(a_n)}{n}=\limsup_{n \to \infty}\dfrac{-v_p(a_n)}{n}$. Therefore, through this paper, we use the symbols $R(\mu):=\{x \in \mathbf{C}_p: v_p(x)>\mu\}$ and $\overline{R}(\mu):=\{x \in \mathbf{C}_p: v_p(x)\geq \mu\}$ to denote region of convergence of power series. Notice that $\mu$ here is considered either finite or $-\infty$ since in the case $\mu=+\infty$, the power series converges only in the point $0$ and there would be nothing to discuss.\\

Let $F(X)=\sum_{i=0}^{\infty}F_iX^i \in \mathbf{C}_p$, we introduce an algorithm to identify a subsequence of indices that possibly determine the number of roots and $p$-adic valuation of the roots.
Constructing the sequence of integer numbers $\{k_0,k_1,\dotsc\}$ with $k_0=0$ and other defined in the manner
\begin{itemize}
	\item Consider the function $l_1(n)=\dfrac{v_p(F_0)-v_p(F_n)}{n}$ with $n \in \{1,2,\dotsc\}$, thus, $\limsup_{n \to \infty}l_1(n)=\mu$. Let $k_1$ be the largest integer such that $l_1(k_1)=\max\{l_1(n): n=1,2,\dotsc \}$. If such the number does not exist which means there are infinitely many integers at that the function $l_1$ obtains the maximal value, then the algorithm ends.
	\item If $\{k_0,k_1,\dotsc,k_{m-1}\}$ is defined then we inductively construct $k_{m}$ which is the largest integer in which the function $l_{m}(n)=\dfrac{v_p(F_{k_{m-1}})-v_p(F_n)}{n-k_{m-1}}$ with $n>k_{m-1}$ is maximal. If such a number does not exist, then the algorithm ends.
\end{itemize}
The algorithm gives a sequence $\{k_0=0,k_1,\dotsc\}$, so-called the characteristic sequence of the power series $F(X)$. In addition, the sequence could be finite or infinite depending on the initial power series. By studying the behavior of this sequence, many properties are revealed.
\begin{theorem}
	Let $F(X)=\sum_{i=0}^{\infty}F_iX^i \in \mathbf{C}_p$ be a power series and its characteristic sequence $\left(k_n\right)_{n \in \textbf{N}}$ induces a strictly decreasing sequence $\left(\alpha_n\right)_{n\geq 1}$ defined by the formula $$\alpha_n=\dfrac{v_p\left(F_{k_{n-1}}\right)-v_p\left(F_{k_{n}}\right)}{k_n-k_{n-1}}$$
	Moreover, the set of roots of $F(X)$ is partitioned into the sets consisting of $k_n-k_{n-1}$ roots of $p$-adic valuation $\alpha_n$.
\end{theorem} 
\begin{proof}
	We use induction on $n$. For $n=1$, we consider the power series
	 $$F\circ p^{\alpha_1}(X)=\sum_{i=0}^{\infty}F_ip^{i\alpha_1}X^i$$
	 By the definition of $k_1$ and $\alpha_1$, we have for any positive integer $n$ 
	 $$\alpha_1 \geq \dfrac{v_p(F_0)-v_p(F_n)}{n} \Leftrightarrow v_p(a_0) \leq v_p(F_n)+n\alpha_1$$
	 and the largest index satisfying the equation is $k_1$, thus, $$\wideg\left(F\circ p^{\alpha_1}(X)\right)=0 \ \text{and} \ \Lwideg\left(F\circ p^{\alpha_1}(X)\right)=k_1$$
	 then the Weierstrass preparation theorem claims that $F\circ p^{\alpha_1}(X)$ has exactly $k_0=0$ roots in open unit disc and $k_1$ roots in closed unit disc. That is interpreted in the manner all $k_1-k_0$ roots of $F\circ p^{\alpha_1}$ are units, thus $F$ has exactly $k_1-k_0$ roots in $\overline{R}(\alpha_1)$, and all of them have $p$-adic valuation $\alpha_1$. Therefore, the statement is correct in the case $n=1$
	 
	 To continue, we use the complete induction on $n$ that is the statement stands for any $n=\{1,\dotsc,m-1\}$ and repeat our reasoning on the function $F\circ p^{\alpha_n}$. Moreover,
	 $$\alpha_{m-1}=\dfrac{v_p\left(F_{k_{m-2}}\right)-v_p\left(F_{k_{m-1}}\right)}{k_{m-1}-k_{m-2}} \ \text{and} \ \alpha_{m-1}>\dfrac{v_p\left(F_{k_{m-2}}\right)-v_p\left(F_{k_{m}}\right)}{k_{m}-k_{m-2}}$$
	 Hence, $$\alpha_{m-1}>\dfrac{v_p\left(F_{k_{m-1}}\right)-v_p\left(F_{k_{m}}\right)}{k_{m}-k_{m-1}}=\alpha_m$$
	 
	 The remain here is that the algorithm discover all roots of the power series $F(X)$. Let $x$ be a root in the region of convergence of $F$ and $v_p(x)=t$, then we consider the power series $F \circ p^t$ has at least a unit root, thus, $\wideg(F \circ p^t)$ and $\Lwideg(F \circ p^t)$ are finite and $t=\dfrac{v_p\left(F_{k_{\wideg(F \circ p^t)}}\right)-v_p\left(F_{\Lwideg(F \circ p^t)}\right)}{k_{\Lwideg(F \circ p^t)}-k_{\wideg(F \circ p^t)}}$ which is also the maximal value of the function $l(\Lwideg(F \circ p^t))$, appearing in the algorithm after at most $\Lwideg(F \circ p^t)$ steps.
\end{proof}
\begin{corollary}
	Every $p$-adic power series has at most countably many roots.
\end{corollary}
Notice that the sequence $\left(\alpha_n\right)_{n >0}$ strictly decreases and admits a lower bound $\mu$., therefore $\lim_{n \to \infty}\alpha_n$ exists.
\begin{proposition}
	$\lim_{n \to \infty}\alpha_n=\mu$
\end{proposition} 
\begin{proof}
	Let indirectly assume $\lim_{n \to \infty}\alpha_n=\eta>\mu$ then by the denseness of $\mathbf{Q}$ in $\mathbf{C}$, there exist a rational number $a$ in the interval $\left(\mu, \eta\right)$, therefore $x$ is a root of $F(X)$ if and only if $p^{-a}x$ is a root of $F\circ p^t(X)$. Furthermore, $F\circ p^t(X)$ converges in $\mathfrak{m_{\mathbf{C}_p}}$ claims that the power series has only finitely many roots in the closed unit disc by the Weierstrass preparation theorem, then so does $F$ over its region of convergence, whence the sequence $\left(\alpha_n\right)_{n >0}$ is finite and $\eta \in \mathbf{Q}$.
	
    Let $n=\wideg(F \circ p^{\eta}(X))$ and $m=\Lwideg(F \circ p^{\eta}(X))$, whence the algorithm stops when processing the function $\dfrac{v_p(F_m)-v_p(F_n)}{n-m}$ with $n >m$. In particular, the function archives maximal value at infinitely many integer, i.e the equality $v_p(F_m)+m\eta=v_p(F_n)+n\eta$ holds for infinitely many integer $n$. However, it implies the region of convergence of $F \circ p^{\eta}(X)$ is $R(0)$ as well as $R(\eta)$ is the region of convergence of $F(X)$ which is contradictory.
\end{proof}
Via what we have accomplished, the condition of a finitely many root $p$-adic power series is determined as the following theorem.
\begin{theorem}
	Let $F(X)$ be power series with coefficients in $\mathbf{C}_p$, then $F(X)$ has infinitely many roots in its region of convergence if and only if the region of convergence is $R(\mu)$ and $\mu$ is either an irrational number or rational number with condition $\liminf_{n \to \infty}v_p(F_n)+n\mu=-\infty$.
\end{theorem}
\begin{proof}
	The power series $F(X)$ has infinitely many roots if and only if the algorithm never ends.
	\begin{itemize}
		\item In the case that $\mu$ is an irrational number then proposition 4.1 ensures that the algorithm would process endlessly.
		\item In the case $\mu$ is a rational number, then there is a correspondence $x \mapsto p^{\mu}x$ between the set of roots of $F(X)$ and $F\circ p^{\mu}(X)$. If the region of convergence of $F\circ p^{\mu}(X)$ is $\mathcal{O}_{\mathbf{C}_p}$ then by the Weierstrass theorem, the number of roots is exactly $\Lwideg(F\circ p^{\mu}(X))$ which is finite. Therefore, $R(\mu)$ is the region of convergence of $F(X)$ and the hypotheses $\wideg(F\circ p^{\mu}(X))$ is finite also means the finiteness of roots, whence  $\wideg(F\circ p^{\mu}(X))=\infty$ as we desire.
	\end{itemize}
\end{proof}
Given a power series $F(X)$ in $\mathbf{C}_p\left[\left[X\right]\right]$ with the characteristic sequence $\{k_0=0,k_1,\dotsc\}$ and the sequence of $p$-adic valuation of roots $\{\alpha_1,\alpha_2,\dotsc\}$, the power series can be factorization in a closed ball $R(t)$ and the circle $C(t)=\{x\in \mathbf{C}_p: v_p(x)=t\}$ with a rational $t$. Let the set of roots $\{z_1,z_2,\dotsc\}$ where for any integer $n$, $v_p(z_i)=\alpha_n \ \forall i \in \{k_{n-1}+1,\dotsc,k_n\}$, let $\mathscr{C}(F,t)=\prod_{\substack{F(z)=0 \\ v_p(z)=t}}(X-z)$ denote the distinguished polynomial in the circle $C(t)$.
\begin{remark}
	The power series $F(X)$ here is considered in the general case, either finitely or infinitely many roots, hence the sequences defined above are finite or infinite respectively. The characteristic above is also known as the Newton polygons of power series which is defined by the terminologies: "slopes", "breaks" and "lengths". Then another interpretation is:
	\begin{itemize}
		\item All breakpoints are of the coordinator $\left(k_i,v_p(F_{k_i})\right)$.
		\item Length of the $i$ slope is $k_i-k_{i-1}$.
		\item The $p$-adic valuations $\alpha_i$ is the inverse slope of Newton polygons.
	\end{itemize}
\end{remark}
\begin{proposition} The polynomial $\mathscr{C}(F,t)$ is determined by coefficients of $F$ for any rational number $t$. Furthermore, if $F(X) \in K(X)$ where $K$ is a finite extension of $\mathbf{Q}_p$ then so is $\mathscr{C}(F,t)$.
	
\end{proposition}
\begin{proof} It is obvious that if $t \notin \{\alpha_1,\alpha_2,\dotsc\}$ then    $\mathscr{C}(F,t)=0$. For any integer $n$, the polynomial $\prod_{i=1}^{k_n}\left(X-\dfrac{z_i}{p^{\alpha_n}}\right)$ is the distinguished polynomial in $\mathcal{O}_{\mathbf{C}_p}$ of the power series $F \circ p^{\alpha_n}(X)$, hence
	$$\prod_{i=1}^{k_n}(X-z_i)=p^{k_n\alpha_n}\mathscr{C}^c\left(F\circ p^{\alpha_n}X\right)\circ p^{-\alpha_n}X$$
    then $$\mathscr{C}(F,\alpha_n)=p^{k_n\alpha_n-k_{n-1}\alpha_{n-1}}\dfrac{\mathscr{C}^c\left(F\circ p^{\alpha_n}X\right)\circ p^{-\alpha_n}X}{\mathscr{C}^c\left(F\circ p^{\alpha_{n-1}}X\right)\circ p^{-\alpha_{n-1}}X}$$
    
    Because the distinguished polynomial in the unit disc is given by the formula as a corollary of the generalization of Washington's theorem, thus so the distinguished polynomial in a circle as well as a closed ball.
    
	Let $z$ be a root of $F(X)$ with $p$-adic valuation $t$, and $Q \in K\left[X\right]$ is the minimal polynomial of $z$, $z'$ is another root of $Q$. By the property of finite extension of $p$-adic valuation, $v_p(z')=v_p(z)$. We observe the $K$-isomorphism $K(z) \to K(z')$ with $z \mapsto z'$ that preserves the $p$-adic valuation, therefore $v_p(F(z))=v_p(F(z'))=\infty$ which also means $z'$ is a root of $F$ with $p$-adic valuation $t$, whence the polynomial $\mathscr{C}(F,t)$ is divisible by $Q$. In fact, the distinguished polynomial is the product of finite minimal polynomials like $Q$, thus $\mathscr{C}(F,t) \in K[X]$
	
\end{proof}
\begin{corollary}
	Roots of any power series with coefficient in a finite extension $K$ of $\mathbf{Q}_p$ are algebraic.
\end{corollary}
\begin{corollary}
	Let $x \in \mathbf{C}_p$ and $\alpha_{n-1} > v_p(x)=t \geq \alpha_n$, then $v_p(F(x))=v_p(F_{k_n})+k_{n-1}t+v_p\left(\mathscr{C}(F,\alpha_n)(x)\right)$
\end{corollary}
\begin{proof}
	Since $v_p(t) \geq \alpha_n$ then $F(t)=F\circ p^{\alpha_n}\left(\dfrac{x}{p^{\alpha_n}}\right)$ and by the theorem 1.2, the $p$-adic valuation equal to 
	$$v_p(F_{k_n})+k_n\alpha_n+v_p\left(\mathscr{C}^c\left(F\circ p^{\alpha_n}\left(\dfrac{x}{p^{\alpha_n}}\right)\right)\right)=\\
	v_p(F_{k_n})+\sum_{i=1}^{k_n}(t-z_i)=v_p(F_{k_n})+k_{n-1}t+v_p\left(\mathscr{C}(F,\alpha_n)(x)\right)$$
\end{proof}
\section{Range of the function $\phi$ in the general case}
In this section, we estimate the maximal value $S$ when $F$ and $G$ are power series in $\mathbf{C}_p\left[\left[X\right]\right]$ over the domains $\mathbf{C}_p$ and the finite extension $K$ of $\mathbf{Q}_p$.
\begin{lemma}
	Let $P(X)=a\prod_{i=1}^{n}(X-z_i)$ and $Q(X)=b\prod_{j=1}^{m}(X-w_j)$ be polynomial in $K\left[X\right]$. The maximal value $S:=\max_{X \in \mathbf{C}_p}\min\left(v_p(P(X)),v_p(Q(X))\right)$ is obtained at one of the points $\{z_1,\dotsc,z_n,w_1,\dotsc,w_m\}$.
\end{lemma}
\begin{proof}
	Let $X$ be an arbitrary point in $\mathbf{C}_p$, because $\overline{\mathbf{Q}_p}$ is dense in $\mathbf{C}_p$ then we can assume $X$ in $K^*$ where $K*$ is a finite extension containing all $z_i$'s, $X$ and $w_j$' with uniformizer $\pi^*$, thus there exists a unique form
	$$X=\sum_{i=-l}^{\infty}x_i\pi^{*i}$$ where $x_i \in \{c_0,c_1,\dotsc,c^{p^f-1}\}$, $x_i$ here is called the $i$-th digit of $X$ in the $\pi^*$-adic expansion. Let $S_k$ be the set given by
	$$S_k:=\{z_i,w_j: z_i \ \text{and} \ w_j \ \text{have the same} \ -l,\dotsc,k \text{-th digits as} \ X \}$$
	It is clear that $S_1 \supset S_2 \supset \dotsc $. This sequence is infinite if and only if $X \in \{z_1,\dotsc,z_n,w_1,\dotsc,w_m\}$
	If $X$ is not a root of $P$ or $Q$ then there exists an integer $k$ such that $S_k \neq \emptyset$ and 
	$$S_{-l} \supset S_{-l+1} \supset \dotsc \supset S_k \supset S_{k+1}= \emptyset $$
	Since $S_k$ is not empty then when we replace $X$ by an element in $S_k$, we assume $z$ then observe the behavior of $v_p(P(X))$ and $v_p(Q(X))$, we have $v_p(P(z)) \geq v_p(Q(X))$ and the same for $Q$. This implies that the value $S$ is obtained at one of the roots of the polynomial $P$ and $Q$.
\end{proof}
\begin{remark}
	The result remains correct if we replace $P$ and $Q$ by two power series since we have already proved that in a closed ball, the $p$-adic valuation of power series only depends on its distinguished polynomial, thus the lemma for power series is true as well as for polynomials.
\end{remark}
\begin{proposition}
	Let $F(X)$ and $G(X)$ be power series with no common root in $\mathbf{C}_p$ and one of them has finitely many roots then the function $\phi$ over the intersection of regions of convergence admits an upper bound.
\end{proposition}
\begin{proof}
	Without loss of generality, we can assume $F$ has exactly $n$ roots in the intersection of regions of convergence and $t$ is the smallest $p$-adic valuation among these roots. Therefore, for any $x$ such that $v_p(x)<t$ then $\phi(x)\leq v_p(F(x))\leq v_p(\Lwideg(F\circ p^t))+\Lwideg(F\circ p^t)t$ which is finite.\\
	For $x \in \overline{R}(t)$, we can the resultant of $P$ and $Q$  in the ball $\overline{R}(t)$ to estimate the maximal value of the function $\phi$ restricted in the closed ball as in section 3.
	 
\end{proof}
Notice that once we can have a method to bound the function $\phi$ restricted to an open ball or restrict ball, the composition shall extend this method to all the intersection of regions of convergence. Therefore, it draws our attention to how to estimate the range of the function $\phi$ in the unit disc(either open or closed).

Let $F$ and $G$ be power series in $\mathbf{C}_p$ of characteristic sequence $\{k_0=0,k_1,\dotsc\}$ and $\{l_0=0,l_1,\dotsc\}$ and the $p$-adic valuation sequences $\{\alpha_1,\alpha_2,\}$ and $\{\beta_1,\beta_2,\dotsc\}$ respectively. The case $F$ and $G$ have common root is obviously eliminated in this paper. We aim to discover when the maximal value is finite and if it is finite, then what is the sharp upper bound of it. It naturally raises a question about the existence of two power series with no common root but $S=+\infty$. The below construction shall answer that matter.
\begin{construction}
	Let consider the sequences of polynomial, for any positive integer $n$
	$$f_n=\prod_{i=1}^{n}p^{-n}\left(X^{2^n}+p\right)$$
	and 
	$$g_n=\prod_{i=1}^{n}\left(p^{2n}+p\right)^{-1}\left(X^{2^n}+p^{2n}+p\right)$$
	Then there exist power series $F$ and $G$ with no common root such that $F$ and $G$ admit all roots of $f_n$ and $g_n$ as roots respectively. In addition, the function $\phi$ is unbounded.
\end{construction}
\begin{proof}
	Notice that for any integer $n$
	$$f_{n+1}=f_n\cdot p^{-1}\left(X^{2^{n+1}}+p\right)=f_n+p^{-1}X^{2^{n+1}}f_n$$
	and 
	$$g_{n+1}=g_n+(p+p^{2n+2})^{-1}X^{2^{n+1}}g_n$$
	Hence, let $f_n \to F$ and $g_n \to G$ in the manner coefficients of $F$ and $G$ are limit of coefficient of $f_n$ and $g_n$ which is already ensured convergent. In addition, the set of roots of $F$ and $G$ are respectively
	$$\bigcup_{n=1}^{\infty}\Bigl\{\text{the}\  n\text{-th roots of}\ -1\ \text{multiplied with} \ p^{1/2^n} \Bigr\}$$ 
	and $$\bigcup_{n=1}^{\infty}\Bigl\{\text{the}\  n\text{-th roots of}\ -1\ \text{multiplied with} \ (p+p^{2n})^{1/2^n} \Bigr\}$$ 
	Hence, $F$ and $G$ do not have any common root. Let $X$ be the $n$-th root of $-1$ multiplied with $p^{1/n}$ and let $n$ tends to $\infty$, then $F(X)=0$, $$v_p(G(X))=v_p(g_n(X))=-n+(n-1)\sum_{i=1}^{n-1}\frac{1}{2^i}+2n \to +\infty$$
	That also means $S=+\infty$.
\end{proof}
\begin{proposition}
	If the intersection if two $p$-adic valuation sequences of $F$ and $G$ is finite then $S$ is finite. i.e $$\left|\{\alpha_1,\alpha_2,\dotsc\} \cap \{\beta_1,\beta_2,\dotsc\}\right|< +\infty \Rightarrow S<+\infty$$
\end{proposition}
\begin{proof}
	Once we consider power series in the unit disc then by proposition 4.1, roots of power series mononically decrease to $0$, therefore $v_p\left(F_{k_{n-1}}\right)>v_p\left(F_{k_{n}}\right)$, then for any $x$ such that $\alpha_{n-1}>v_p(x)>\alpha_n$, then $\phi(x)< v_p(F_{k_n})+k_{n-1}v_p(x)+(k_n-k_{n-1})\alpha_n$ by corollary  4.3 which is bounded by 
	$$v_p(F_{k_n})+k_{n-1}v_p(x)+(k_n-k_{n-1})\alpha_n<v_p(F_{k_{n-1}})+k_{n-1}\alpha_{n-1} \leq v_P(F_0)$$
	Combining with the lemma 5.1, we obtain 
	$$S \leq \max\{v_p(F(0)),v_p(G(0)),\phi(x): x \in \{\alpha_1,\alpha_2,\dotsc\} \cap \{\beta_1,\beta_2,\dotsc\} \}<+\infty$$
\end{proof}
To continue our estimation, we observe the function $\phi$ in the circle $C(\alpha_n=\beta_m)$ for some indices $n$ and $m$.
\begin{theorem}
	Let $n$ and $m$ be integer such that $\alpha_n=\beta_m$, $a=\min\{v_p(F_{k_n})+k_{n-1}\alpha_n,v_p(F_{l_m})+l_{m-1}\alpha_n \}$, $b=\max\{v_p(F_{k_n})+k_{n-1}\alpha_n,v_p(F_{l_m})+l_{m-1}\alpha_n \}$ and $r=\Res(\mathscr{C}(F,\alpha_n),\mathscr{C}(G,\beta_m))$, then 
	$$a+\dfrac{v_p(r)}{\min\{l_m-l_{m-1},k_n-k_{n-1}\}}\leq\max_{x \in C(\alpha_n)}\phi(x)\leq b+ v_p(r)$$
	Furthermore, if $F$ and $G$ are in $K\left[\left[X\right]\right]$ where $K$ is a finite extension with the ramification index $e$ and uniformizer $\pi$, then the upper bound can be sharpened as following 
	$$\max_{x \in C(\alpha_n)}\phi(x)\leq  b+cv_p(r)$$
	where $c=\max\Bigl\{\dfrac{\gcd(v_{\pi}(F_{k_n})-v_{\pi}(F_{k_{n-1}}),k_n-k_{n-1})}{k_n-k_{n-1}},\dfrac{\gcd(v_{\pi}(G_{l_m})-v_{\pi}(G_{l_{m-1}}),l_m-l_{m-1})}{l_n-l_{m-1}}\Bigr\}$
\end{theorem}
\begin{proof}
	By using corollary 4.3 and lemma 5.1, it is enough to prove that let $\{z_i: i=\overline{k_{n-1}+1,k_n}\}$ and $\{w_j: j=\overline{l_{m-1}+1,l_m}\}$ be the set of $\mathscr{C}(F,\alpha_n)$ and $\mathscr{C}(G,\beta_m)$, then 
	$$\dfrac{v_p(r)}{k_n-k_{n-1}}\leq \max\{v_p(\mathscr{C}(F,\alpha_n)(w_j)): j=\overline{l_{m-1}+1,l_m}\}\leq v_p(r) \ (1)$$
	and $$\dfrac{v_p(r)}{l_m-l_{m-1}}\leq \max\{v_p(\mathscr{C}(F,\alpha_n)(z_i)): i=\overline{k_{n-1}+1,k_n}\}\leq v_p(r) \ (2)$$
	In the case $F$ and $G$ are in $K\left[\left[X\right]\right]$,
	$$\max\{v_p(\mathscr{C}(F,\alpha_n)(w_j)): j=\overline{l_{m-1}+1,l_m}\}\leq \dfrac{\gcd(v_{\pi}(F_{k_n})-v_{\pi}(F_{k_{n-1}}),k_n-k_{n-1})}{k_n-k_{n-1}}v_p(r) \ (3)$$
	and $$\max\{v_p(\mathscr{C}(F,\alpha_n)(z_i)): i=\overline{k_{n-1}+1,k_n}\}\leq \dfrac{\gcd(v_{\pi}(G_{l_m})-v_{\pi}(G_{l_{m-1}}),l_m-l_{m-1})}{l_n-l_{m-1}}v_p(r) \ (4)$$
	The inequalities $(1)$ and $(2)$ are obtained directly from the definition of resultant
	$$v_p(r)=v_p\left(\prod_{j=l_{m-1}+1}^{l_m}\mathscr{C}(F,\alpha_n)(w_j)\right)=
	v_p\left(\prod_{i=k_{n-1}+1}^{k_n}\mathscr{C}(G,\beta_m)(z_i)\right)$$
	Using the condition $F \in K\left[\left[X\right]\right]$, then $\mathscr{C}(F,\alpha_n)(X) \in \mathcal{O}_K\left[X\right]$ which is a monic polynomial of degree $k_n-k_{n-1}$ and has the constant coefficient of $p$-adic valuation $v_{\pi}(F_{k_n})-v_p(F_{k_{n-1}})$, therefore every root of $\mathscr{C}(F,\alpha_n)(X)$ must have the degree of its minimal polynomial at least $\dfrac{k_n-k_{n-1}}{\gcd(v_{\pi}(F_{k_n})-v_{\pi}(F_{k_{n-1}}),k_n-k_{n-1})}$, whence we achieve $(3)$ and $(4)$.
\end{proof}
\begin{proposition}
	If one of the following sets is bounded then $S < +\infty$:
	\begin{itemize}
		\item $\Bigl\{v_p\left(\Res(\mathscr{C}(F,\alpha_n),\mathscr{C}(G,\beta_m))\right): \alpha_n=\beta_m\Bigr\}$
		\item $\Bigl\{\dfrac{\gcd(v_{\pi}(F_{k_n})-v_{\pi}(F_{k_{n-1}}),k_n-k_{n-1})}{k_n-k_{n-1}}v_p\left(\Res(\mathscr{C}(F,\alpha_n),\mathscr{C}(G,\beta_m))\right),\\
		\dfrac{\gcd(v_{\pi}(G_{l_m})-v_{\pi}(G_{l_{m-1}}),l_m-l_{m-1})}{l_n-l_{m-1}}v_p\left(\Res(\mathscr{C}(F,\alpha_n),\mathscr{C}(G,\beta_m))\right): \alpha_n=\beta_m\Bigr\}$
	\end{itemize}
		Furthermore, if there exist a upper bound of the set $$\Bigl\{\gcd(v_{\pi}(F_{k_n})-v_{\pi}(F_{k_{n-1}}),k_n-k_{n-1}),\gcd(v_{\pi}(G_{l_m})-v_{\pi}(G_{l_{m-1}}),l_m-l_{m-1}): \alpha_n=\beta_m\Bigr\}$$ then $S$ is finite if and only if $$\sup\Bigl\{\dfrac{v_p\left(\Res(\mathscr{C}(F,\alpha_n),\mathscr{C}(G,\beta_m))\right)}{k_n-k_{n-1}},\dfrac{v_p\left(\Res(\mathscr{C}(F,\alpha_n),\mathscr{C}(G,\beta_m))\right)}{l_m-l_{m-1}}: \alpha_n=\beta_m\Bigr\}< +\infty$$
	
\end{proposition}
\begin{proof}
	These statements are obtained directly via the above inequalities. 
\end{proof}
\section{Applications of resultant to irreducibility and roots of polynomials  }
In this section, we utilize the expansion of elements of $\overline{\mathbf{Q}}_p$ to decide the irreducibility and more generally to approximate roots of given power series. First, we mention known irreducibility criteria in terms of what we have achieved in this paper.   
\begin{proposition}
	Let $P(X)=X^n+a_{n-1}X^{n-1}+\dotsc+a_0$ be a polynomial with coefficients in an extension $K$ and the characteristic sequence $\{k_0=0,k_1,\dotsc,k_m\}$ then:
	\begin{itemize}
		\item $m>1$ then $P$ is reducible.
		\item If $m=1$ and $v_{\pi}\left(a_0\right)$ and $n$ are relatively prime then $P$ is irreducible over $K$.
	\end{itemize}	
\end{proposition} 
\begin{proof}
	The first two statements are obtained directly from the properties of the characteristic sequence. In terms of Newton polygons, those are interpreted as $P$ is irreducible if the Newton polygon is pure (i.e consisting of only one slope) and there is no point in that slope which is a stronger version of Einstein criterion.
\end{proof}

Let $K$ be a finite extension then we use notations $e_K$, $f_K$ and $\mathrm{k}_K$ to denote the ramification index, residue degree and residue field respectively. Through this section, by a set of representative of the residue field, we mean the finite field $\mathrm{F}_{p^f}$ where $f$ is the residue degree. Let $\overline{\mathbf{Q}}_p^*$ be the set of all algebraic over $\mathbf{Q}_p$ having the degree of minimal polynomial non-divisible by $p$ and the map $*:\overline{\mathbf{Q}}_p^* \to\overline{\mathrm{F}}_p \times\mathbf{Q}$ defined by 
$$ x \mapsto \varepsilon p^t$$ where $v_p(x)=t \in \mathbf{Q}$ and $\varepsilon \in \overline{\mathrm{F}}_p$ such that $v_p(x-x^*) > t$. The existence of this map is provided by the following lemmas.

\begin{lemma}
	Let $P(X)=X^n+a_{n-1}X^{n-1}+\dotsc+a_0$ be a polynomial over $\mathbf{Q}_p$ such that $a_0$ is a unit and the factorization of $P$ over $\mathrm{F}_p$ is $P_1\dotsc P_m$, then the followings are obtained:
	\begin{itemize}
		\item[i)] If $P$ is irreducible over $\mathrm{F}_p$ then $P$ is irreducible over $\mathbf{Q}_p$.
		\item[ii)] For any integer $i=\overline{1,m}$, $P$ has a root $x$ satisfying $f_{\mathbf{Q}_p(x)}$ divisible by $\deg(P_i)$.
		\item[iii)] If $P$ is irreducible then $P_1=\dotsc=P_n$, i.e the polynomial $P$ over $\mathrm{F}_p$ is a power of an irreducible $Q(X) \in \mathrm{F}_p$.
	\end{itemize}
\end{lemma}
\begin{proof}
	Notice that if $P$ is reducible over $\mathbf{Q}_p$ then consequently that is reducible over $\mathrm{F}_p$, and the irreducibility of a polynomial over a finite field is possibly checked by Rabin's test $\cite{PA}$. By factorizing the polynomial $P$ over the algebraic closure $\overline{\mathbf{Q}}_p$ 
	$$P(X)=\prod_{i=1}^{n}(X-x_i)$$ then it is clear that the polynomial $P$ over admits the below factorization over the residue field $\overline{\mathrm{F}}_p$
	$$P(X)=\prod_{i=1}^{n}(X-x_i')$$  
	where $x_i' \in \overline{\mathrm{F}}_p$ which is also in $\mathbf{Q}_p(x_i)$, whence statement $ii$ is attained.
	Assuming that $P(X)$ is irreducible, then for integers $i\neq j \in \{1,\dotsc,n\}$, the $\mathbf{Q}_p$-homomorphism $\mathbf{Q}_p(x_i) \to \mathbf{Q}_p(x_j)$ defined by $x_i \to x_j$ preserves the $p$-adic valuation. Without loss of generality, let $x_i'$ and $x_j'$ be roots of $P_1$ and $P_2$ respectively. Therefore, $P_1(x_i)$ is not a unit in $\mathbf{Q}_p(x_i)$ which implies the same about $P_1(x_j)$ in $\mathbf{Q}_p(x_j)$ , then $P_1(x_j')=0$. However, $P_1$ and $P_2$ are irreducible over $\mathrm{F}_p$, whence $P_1=P_2$. By the same reasoning, we have $P_1=\dotsc=P_m=Q$ and $P(X)\equiv(Q(X))^m$ mod $(p)$.

\end{proof}
\begin{corollary}
	It is obvious that for any root $x$ of $P(X)$, $\deg(Q)\mid f_{\mathbf{Q}_p(x)}$ and $e_{\mathbf{Q}_p(x)} \mid m$ those imply a new irreducibility criterion that a monic polynomial with coefficients in $\mathbf{Q}_p$ and constant one is a unit, if the polynomial is not a power mod $(p)$, then it is reducible.
\end{corollary}
\begin{lemma}
	Let $p \nmid n$ and $P(X)=X^n+a_{n-1}X^{n-1}+\dotsc+a_0$ be a polynomial over $\mathbf{Q}_p$ and every root of $P$ has the same $p$-adic valuation $t=\dfrac{v_{p}(a_0)}{n}>0$ and the factorization of $P$ over the quotient ring $\mathcal{O}_K/(p^{v_{p}(a_0)+1})$ is $P_1\dotsc P_m$, then the followings are obtained:
	\begin{itemize}
		\item[i)] If $P$ is irreducible over the quotient ring $\mathcal{O}_{\mathbf{Q}_p}/(p^{v_{p}(a_0)+1})$ then $P$ is irreducible over $\mathbf{Q}_p$.
		\item[ii)] If $P$ is irreducible then $P$ is a power of a irreducible polynomial over $\mathcal{O}_{\mathbf{Q}_p}/(p^{v_{p}(a_0)+1})$.
	\end{itemize}
\end{lemma}
\begin{proof}
The statement $i$ is trivial and let $t=u/v$ where $u$ and $v$ are co-prime natural numbers.

If the polynomial $P$ is irreducible then for any $x$ and $x'$ roots, we have a $\mathbf{Q}_p$-isomorphism $$\mathbf{Q}_p(x) \to \mathbf{Q}_p(x')$$ $$x \mapsto x'$$
Besides, if we take any element $y$ of $\overline{\mathbf{Q}}_p$ such that $v_p(x-y) >t$ then $P(y)$ vanishes mod $(p^{v_p(a_0)+1})$ as well. Therefore, we can consider the polynomial $Q=\dfrac{P\circ p^tX}{p^{tn}}$ embedded into $\mathcal{O}_{\mathbf{Q}_p}^{\times}$ which is a monic and all roots are units. By that, we mean omitting all non-unit coefficients of $Q$ and replacing the remains by their last digit which is the congruent mod $(p)$ root of unity. Let the last digits of all roots be $\varepsilon_i$, $i=\overline{1,n}$ then we have those coefficients with degree non-divisible by $u$ are omitted in the process. Hence, there exists a polynomial $Q'$ of degree $n/v$ such that $Q(X)=Q'(X^v)$ and all roots of $Q'$ are units, let $\varepsilon_i' \ i=\overline{1,n/v}$ be the last digits of those roots then $$Q'(X)=\prod_{i=1}^{n/v}(X-\varepsilon_i')+pH(X), \ H(X)\in \mathcal{O}_{\mathbf{Q}_p}[X] $$
Moreover, since $p \nmid n$ then there exist roots $\varepsilon_i, \ i=\overline{1,n}$ of unity such that $$Q(X)=\prod_{i=1}^{n}(X-\varepsilon_i) +pH(X^v)$$
It is clear that $\varepsilon_i$'s are the last digits of roots of the polynomial $P$. Furthermore,$P(X)=\prod_{i=1}^{n}(X-\varepsilon_i)+pH'(X), \ H'(X)\in \mathcal{O}_{\mathbf{Q}_p}[X]$, i.e for any root $x$ of $P$ there exists an index $i$ such that $v_p(x-\varepsilon_ip^t)>t$ which combines with lemma 6.1 implies the map $*$ is well-defined. Now, let the integer $i$ such that $v_p(x-\varepsilon_ip^t)>t$, then we extend to the $\mathbf{Q}_p$-homomorphism 
$$\begin{tikzcd}
	\mathbf{Q}_p(x) \arrow{r}{x \mapsto x'} \arrow[swap]{d}{} & \mathbf{Q}_p(x',\varepsilon_1p^t,\dotsc,\varepsilon_np^t) \arrow{d}{} \\%
	\mathbf{Q}_p(x,\varepsilon_ip^t) \arrow{r}{\varphi}& \mathbf{Q}_p(x',\varepsilon_1p^t,\dotsc,\varepsilon_np^t)
\end{tikzcd}$$
Since $v_p(x-\varepsilon_ip^t) > t \Rightarrow v_p(x'-\varphi(\varepsilon_ip^t)) > t$ and $\varphi(\varepsilon_ip^t) \in \{\varepsilon_1p^t,\dotsc,\varepsilon_np^t\}$, whence $\varphi(\varepsilon_ip^t)=x'^*$ and $\mathbf{Q}_p(x,x^*) \simeq \mathbf{Q}_p(x',x'^*)$ by an isomorphism which maps $x^*$ to $x'^*$. The polynomial $\prod_{i=1}^{n}(X-\varepsilon_i p^t) \in \mathbf{Q}_p[X]$ has an isomorphism between any two roots, then its is a power of irreducible polynomial. 
\end{proof}
\begin{theorem}
	Let $x$ be an algebraic of the field $\mathbf{Q}_p$ of non-negative $p$-adic valuation and $P(X)=X^n+a_{n-1}X^{n-1}+\dotsc+a_0 \in \mathbf{Q}_p[X]$ be the minimal polynomial of $x$ of degree $p \nmid n$, then there exist a monotonically increasing integer sequence $\{t_i=u_i/v_i: u_i, v_i \in \textbf{Z},\ \gcd(u_i,v_i)=1\}$ and the sequence $\{\varepsilon_i\ \text{is an algebraic of}\ \mathrm{F}_p , i \in \textbf{Z}\}$ in $\mathrm{F}_p$ such that $$x:=\sum_{i=0}^{\infty}\varepsilon_ip^{t_i}$$
	and the followings are satisfied
	
	\begin{itemize}
		\item[i)] $\lim_{i \to \infty}\lcm(v_1,\dotsc,v_i)=e_{\mathbf{Q}_p(x)}$ and $\varepsilon_i^{v_i} \in \mathbf{Q}_p(x)$.
		\item[ii)] There is a finite extension of $\mathrm{F}_p$ containing all $\varepsilon_i$.
		\item[iii)] The expansion is unique.
    \end{itemize}
\end{theorem}
\begin{proof}
	Let $P(X)=\prod_{i=1}^{n}(X-x_i)$ with $x=x_1$ and run an algorithm as following on the roots of $P(X)$.
	\begin{itemize}
		\item If $a_0$ is a unit, we consider the polynomial $P^{(0)}=P \mod p$ over the residue field $\mathrm{k}_{\mathbf{Q}_p}$ and its factorization. It is clear that the set of last digits of all roots of $P$ are all roots of $P^{(0)}$. Let $\varepsilon_0$ be a root of unity and the last digit of $x$, i.e $x-\varepsilon_0$ is not a unit. If $a_0$ is not a unit then $\varepsilon_0=0$.
		\item Translate the polynomial $P$ by $-\varepsilon_0$ which is $P(X+\varepsilon_0)$, let $t_1$ be $p$-adic valuation of $x-\varepsilon_0$ which is a root of $P(X+\varepsilon_0)$. Hence, the distinguished polynomial in the circle $C(t_1)$ is a polynomial with coefficients in $\mathbf{Q}_p(\varepsilon_0)$ denoted by $P^{(1)}$ having all roots of $p$-adic valuation $t_1$.Therefore, $(x-\epsilon_0)^*$ is a root of $P^{(1)}$ mod $(p)$. Then $\epsilon_1p^{t_1}=(x-\epsilon_0)^*$ by lemma 6.2.
		\item The polynomial $P(X+\varepsilon_0+\varepsilon_1)$ over the field $\mathbf{Q}_p(\varepsilon_0,\varepsilon_1p^{t_1})$ to obtain $\varepsilon_2p^{t_2}=(x-\varepsilon_0-\varepsilon_1)^*$.
		\item Repeat the process.
	\end{itemize}
	We observe the algorithm. The lemma 6.1 claims that $P^{(0)}(X)=(Q^{(0)})^{m_0}$ over the residue field $\mathrm{k}_{\mathbf{Q}_p}$, then $m_0 \mid n$. In the case $a_0$ is not a unit then $m_0=n$ and $Q^{(0)}\equiv 0$. Additionally, $Q^{(0)}$ is the minimal polynomial of $\varepsilon_0$ over $\mathrm{F}_p$ and there are exactly $m_0$ roots of $P$ of the last digit $\varepsilon_0$. Moreover, for the roots $x_i$ and $x_j$, we have $\mathbf{Q}_p$-isomorphism $$\mathbf{Q}_p(x_i) \to \mathbf{Q}_p(x_j)$$
	This map preserves $p$-adic valuation, whence $v_p(x_i-\varepsilon_0)=t_1 \forall i$ such that $v_p\left(x_i-\varepsilon_0\right)>0$. That means $\deg(P^{(1)})=m_0 \mid n$.
	
	Let $m_l$ be the number of $x_j$ such that $$v_p\left(x_j-\sum_{i=0}^{l}p^{t_i}\varepsilon_i\right)>t_l$$
	\begin{lemma}
		The polynomial $P^{(l)}$ defined in the $l$-th step of the algorithm is a $m_l$ degree power of an irreducible polynomial $Q^{(l)}$ over $\mathrm{F}_p(\varepsilon_0,\dotsc,\varepsilon_{l-1})$. Furthermore, $m_l$ is also the degree of $P^{(l+1)}$.
	\end{lemma}
	\begin{proof}
		We use induction on $l$ and the case $l=0$ is already proved above. Let assume the statement is correct for all $0,\dotsc,l$.  Let $x_i$ and $x_j$ be any two of those $m_l$ roots, then one can be extended to a $\mathbf{Q}_p(\varepsilon_0,\dotsc,\varepsilon_l,p^{t_1},\dotsc,p^{t_l})$-isomorphism
		$$\begin{tikzcd}
			\mathbf{Q}_p(x_i) \arrow{r}{\simeq} \arrow[swap]{d}{} & \mathbf{Q}_p(x_j) \arrow{d}{} \\%
			\mathbf{Q}_p(x_i,\varepsilon_0,\dotsc,\varepsilon_lp^{t_l}) \arrow{r}{\simeq}& \mathbf{Q}_p(x_j,\varepsilon_0,\dotsc,\varepsilon_lp^{t_l})
		\end{tikzcd}$$
	Hence, $v_p\left(x_i-\sum_{i=0}^{l}p^{t_i}\varepsilon_i\right)=v_p\left(x_j-\sum_{i=0}^{l}p^{t_i}\varepsilon_i\right)=t_{l+1}\Rightarrow\left(x_i-\sum_{i=0}^{l}p^{t_i}\varepsilon_i\right)^*\mapsto \left(x_j-\sum_{i=0}^{l}p^{t_i}\varepsilon_i\right)^*$. That implies $P^{(l+1)}$ is a power of degree $m_{l+1}$ which is the degree of $P^{(l+2)}$ as well.
	\end{proof}
	The lemma provide a sequence of integer $(m_l)_{l \geq 0}$ and $m_{l+1} \mid m_l \forall l$ which shows that the algorithm shall end after finite steps because the sequence $(m_l)_{l \geq 0}$ stabilizes from some index, i.e under finite translations, the polynomial 
	$$P^{(l+1)}\left(X+\sum_{i=0}^{l}p^{t_i}\varepsilon_i\right)=\left(X-\sum_{j=l+1}^{\infty}p^{t_j}\varepsilon_j\right)$$
	for some $l$.

    The algorithm gives a result $$x:=\sum_{i=1}^{\infty}\varepsilon_ip^{t_i}$$
    Besides, we already know that $v_p\left(\mathbf{Q}_p(x)\right) \subset \frac{1}{n}\textbf{Z}$, and we use introduction on index to prove that $\lcm\left(v_1,\dotsc,v_i\right) \mid e_{{\mathbf{Q}_p(x)}}$.
    
    It is clear that $\varepsilon_0 \in \mathrm{k}_{\mathbf{Q}_p(x)} \in \mathbf{Q}_p(x)$, thus, $x-\varepsilon_0 \in \mathbf{Q}_p(x)$, whence $t_1 \in  v_p\left(\mathbf{Q}_p(x)\right)$ which also means $v_1 \mid e_{{\mathbf{Q}_p(x)}}$. Notice that 
    $$\left(x-\varepsilon_0\right)^{v_1}=\varepsilon_1^{v_1}p^{u_1}+v_1(\varepsilon_1p^{t_1})^{v_1-1}\varepsilon_2p^{t_2} +\text{multiple of}\ p^{u_1+2t_2-2t_1})$$	
    Dividing with $p^{u_1}$, we get $\varepsilon_1^{v_1}$ and $v_2 \mid e_{{\mathbf{Q}_p(x)}}$. By induction, we get $\lim_{i \to \infty}\lcm(v_1,\dotsc,v_i)\mid e_{\mathbf{Q}_p(x)}$
    Let $c=\lim_{i \to \infty}\lcm(v_1,\dotsc,v_i)$ and notice that any element $y$ of $\mathbf{Q}_p(x)$ takes the form $$y:= \dfrac{f(x)}{g(x)}$$ where $f$ and $g$ are polynomials in $\mathbf{Q}_p[X]$, then $v_p(y)$ is multiple of $c$, that means $e_{\mathbf{Q}_p(x)} \mid c$, whence the statement $i$ and $ii$ are obtained.
    
    Combining $i$ and $ii$, we have $\mathbf{Q}_p(x) \in \mathbf{Q}_p(p^{1/n},\varepsilon)$ where $\varepsilon$ is the $p^{n\varphi(n)}-1$-th primitive root of unity since $\varepsilon_i^n \in \mathrm{k}_{\mathbf{Q}_p(x)}$. Therefore, the statement $iii$ is proved because $p^{1/n}$ is a uniformizer of $\mathbf{Q}_p(p^{1/n},\varepsilon)$.
    Furthermore, the sequences $(t_i)$ of  Galois conjugates $x'$ of $x$ is the same because of extension   $$\begin{tikzcd}
    	\mathbf{Q}_p(x) \arrow{r}{\simeq} \arrow[swap]{d}{} & \mathbf{Q}_p(x') \arrow{d}{} \\%
    	\mathbf{Q}_p(x,\varepsilon_0,\dotsc,\varepsilon_lp^{t_l}) \arrow{r}{\simeq}& \mathbf{Q}_p(x',\varepsilon_0',\dotsc,\varepsilon_l'p^{t_l})
    \end{tikzcd}$$
    
\end{proof}
\begin{corollary}
	Let $K$ be a finite extension and $[K:\mathbf{Q}_p]=m$ non-divisible by $p$, then $K \subset\mathbf{Q}_p(p^{1/m},\varepsilon)$ where $\varepsilon$ is the $p^{m}-1$-th primitive root of unity.
\end{corollary}
\begin{corollary}
	The above algorithm raises a question about how many steps the translated polynomial $P^{(m)}$ is a linear polynomial. The resultant gives th upper bound $nv_p\left(\Res(P,P')\right)$.
\end{corollary}
\begin{proof}
	Let $x:=\sum_{i=1}^{\infty}\varepsilon_ip^{t_i}$ be a fixed roots,
	because $P$ is irreducible then all roots are distinct. That means there is the smallest integer $m$ such that $(\varepsilon_i: i=0,\dotsc,m) \neq  (\varepsilon_i': i=0,\dotsc,m) \forall x':=\sum_{i=1}^{\infty}\varepsilon_i'p^{t_i}$ but there exist a root $x'':=\sum_{i=1}^{\infty}\varepsilon_i''p^{t_i}$ satisfying $(\varepsilon_i: i=0,\dotsc,m-1) \neq  (\varepsilon_i'': i=0,\dotsc,m-1)$. Consequently, $P^{(m)}$ is a linear polynomial. Besides, since $t_i >1/n$, then $m/n<v_p(x-x'')<v_p\left(\Res(P,P')\right) \Rightarrow m <nv_p\left(\Res(P,P')\right)$.
\end{proof}

By putting all pieces together, we develop algorithms to determine roots and the irreducibility of a given polynomial as the following.
\begin{algorithm}
	The algorithm in proof of the theorem 6.1 can be corrected to find roots of an arbitrary polynomial $P$ in the manner 
	\begin{itemize}
	\item Breaking down $P$ into a product of polynomial of the same $p$-adic valuation root.
	\item For each obtained polynomial, continue dividing it into the spare-free polynomial.
	\item For each obtained polynomial above, finding digits of roots by the algorithm in proof of theorem 6.1.
	\end{itemize}
\end{algorithm}
\begin{algorithm}
	Given a polynomial $P(X)$ of degree $n$ non-divisible by $p$ with coefficient in $\mathbf{Q}_p$, a test of irreducibility can be developed from what we already have.
		\begin{itemize}
		\item Check the polynomial is square-free or not by resultant.
		\item Consider the square-free polynomial $P(X)$ and use the algorithm in the proof of theorem 6.1.
	\end{itemize}
In each step, we use Berlekamp's algorithm $\cite{BR}$ to the polynomials $F^{l}$ as above to define if it is a power over a finite field. If there is a step the condition is not satisfied then the polynomial is reducible. If not, after finite steps, we shall get the translated polynomial 
	$$P^{(l+1)}\left(X+\sum_{i=0}^{l}p^{t_i}\varepsilon_i\right)=\left(X-\sum_{j=l+1}^{\infty}p^{t_j}\varepsilon_j\right)^m$$
	If $m=1$ then $P$ is irreducible and if $m>1$ then $P$ is a $m$ degree power of an irreducible polynomial, whence it is reducible. As above, the number of steps can be bounded by $nv_p\left(\Res(P,P')\right)$.
\end{algorithm}
\section*{Acknowledgement}
I wish to thank my supervisor Dr. Zábrádi Gergely for his assistance during the time I study the topic.

\newcommand{\Addresses}{}
\bigskip
\footnotesize
\textsc{Department of Mathematics,\
	Eötvös Loránd University}\par\nopagebreak
\textit{E-mail address}, Tran Hoang Anh: \texttt{bongtran5399@gmail.com}

\end{document}